\documentclass[a4paper]{amsart}

\usepackage{amssymb,amsmath,latexsym,epsfig,euscript, amsthm}

\def\boxit#1#2{\setbox1=\hbox{\kern#1{#2}\kern#1}%
\dimen1=\ht1 \advance\dimen1 by #1 \dimen2=\dp1 \advance\dimen2 by
#1
\setbox1=\hbox{\vrule height\dimen1 depth\dimen2\box1\vrule}%
\setbox1=\vbox{\hrule\box1\hrule}%
\advance\dimen1 by .4pt \ht1=\dimen1 \advance\dimen2 by .4pt
\dp1=\dimen2 \box1\relax}

\def\R{\mathbb{R}}
\def\Q{\mathbb{Q}}
\def\Z{\mathbb{Z}}
\def\C{\mathbb{C}}
\def\A{\mathbb{A}}
\def\N{\mathbb{N}}

\def\T{\mathbb{T}}

\def\OO{\mathcal{O}}

\def\div{\, |\,}
\def\ndiv{\nmid}

\def\ess{{\rm ess}}

\def\Qbarra {{\overline{\Q}}}
\theoremstyle{definition}
\newtheorem{defn}{Definition}[section]

\theoremstyle{plain}
\newtheorem{lem}[defn]{Lemma}

\newtheorem{prop}[defn]{Proposition}
\newtheorem{teo}[defn]{Theorem}

\newtheorem{cor}[defn]{Corollary}

\newtheorem{conjecture}[defn]{Conjecture}

\newenvironment{undef}[1]%
           {\vspace{3.3mm}
           \noindent{\bf #1}\it}%
           {\vspace{3.3mm}}

\newenvironment{Proof}[1]{
  \trivlist \item[\hskip \labelsep{\it #1.}]}{\hfill\mbox{$\Box$}
  \endtrivlist}

\begin{document}

\title{Factoring bivariate sparse (lacunary) polynomials}

\author{Mart\'\i n Avenda\~no }
\address{Departamento de Matem\'atica, Facultad de Ciencias Exactas y
Naturales, Universidad de Buenos Aires}
\email{mavendar@bigua.dm.uba.ar}
\thanks{M. Avenda\~no was  supported by a CONICET fellowship, Argentina.}
\author{Teresa Krick}
\address{Departamento de Matem\'atica, Facultad de Ciencias Exactas y
Naturales, Universidad de Buenos Aires} \email{krick@dm.uba.ar}
\thanks{T. Krick was partially supported by research grants UBACYT
X-112 and CONICET PIP 2461/01, Argentina.}
\author{Mart\'\i n Sombra}
\address{Departament d'\`Algebra i Geometria, Universitat de Barcelona} \email{sombra@ub.edu}
\thanks{M. Sombra was supported by the Ram\'on y Cajal program of the Ministerio de Educaci\'on y Ciencia,  Spain.}

\date{\today}

\begin{abstract} We present a deterministic algorithm for computing
all irreducible factors of degree $\le d$ of a given bivariate
polynomial $f\in K[x,y]$ over an algebraic number field $K$ and
their multiplicities, whose running time is polynomial in the bit
length of the sparse encoding of the input and in $d$. Moreover, we
show that the factors over $\Qbarra$ of degree $\le d$ which are not
binomials  can also be computed in time polynomial in the sparse
length of the input  and in  $d$.
\end{abstract}
\maketitle


\typeout{Introduccion}
\section*{Introduction}

Effective factorization of polynomials, when possible,  is an
important task in computational algebra and number theory. This
problem has a long story, going back to I. Newton in 1707, and to
the astronomer F. von Schubert who in 1793 presented an algorithm
for factoring a univariate polynomial, later rediscovered and
generalized by  L. Kronecker in 1882. Many other more efficient
algorithms were designed since then, between the most famous ones we
cite~\cite{Ber70,Zas69}.

In 1982, A.K. Lenstra, H.W. Lenstra Jr. and L. Lov\'asz made a
fundamental advance by obtaining the first deterministic
polynomial-time algorithm for factoring a univariate polynomial over
the rationals~\cite{LLL82}. Based on this result and the technique
of lattice basis reduction introduced for its proof, several new
factorization algorithms were
obtained~\cite{ChGr82,Len84,Kal85,Lan85,Len87}. These algorithms
succeeded in bringing to polynomial time the problem of factoring
univariate and multivariate polynomials over algebraic number fields
when given by their {\it dense} encoding, that is the input $f$ is
given by the list of all its terms of degree $\le \deg (f)$
including the zero ones.

\medskip
For practical purposes, it is usually more realistic to consider the
{\it sparse} (or lacunary) encoding of a polynomial. In this paper
we consider the problem of factoring a {bivariate} polynomial
$$
f=\sum_{i=1}^t a_ix^{\alpha_i}y^{\beta_i} \in \Q[x,y]
$$
given in {sparse}  encoding: by the list
$(a_i,\alpha_i,\beta_i)_{1\le i\le t}$ of its non-zero coefficients
and corresponding exponents. Let $\ell(f)$ denote the {\it bit
length} of the sparse encoding of $f$; informally speaking this is
the number of bits needed to spell out the data. We obtain a
deterministic algorithm for computing the low degree factors of $f$
in time polynomial in $\ell(f)$:

\begin{undef}{Theorem 1.}
There is a deterministic algorithm that, given $f\in \Z[x,y]$ and $d
\ge 1$, computes all irreducible factors of $f$ in $\Q[x,y]$ of
degree $\le d$ together with their multiplicities, in
$(d\cdot\ell(f))^{O(1)}$ bit operations.
\end{undef}

More generally, this algorithm applies for factoring bivariate
polynomials over number fields~(see Subsection~\ref{Rational
factorisation}).

\medskip
Let us observe that since the degree of a polynomial can be
exponentially big in its sparse length (we have $ \deg(f) \le
2^{\ell(f)} $ and this upper bound is attainable),  a direct
application of the algorithms for factoring dense polynomials would
give an exponential complexity. The restriction to bounded degree
factors is unavoidable: the polynomial $f=x^p-1$ ($p$ prime) is of
sparse length $\log_2 (p) + O(1) $ but has the  dense irreducible
factor $x^{p-1}+\cdots + 1$.

\medskip   The first result in this
direction  appeared in 1998, when F. Cucker, P. Koiran and S. Smale
showed how to find all the integer roots of a univariate polynomial
with integer coefficients in  polynomial time in its sparse
encoding, and asked whether one can find in the same time the
rational roots as well~\cite{CKS99}. This question (and more!) was
affirmatively answered by  H.W. Lenstra Jr. who presented an
algorithm that
---given a number field $K$ and a univariate polynomial $f\in K[x]$--- computes
all its irreducible factors of degree $\le d$ together with their
multiplicities, in $(d+\ell(f))^{O(1)}$ bit
operations~\cite[Thm]{Len99b}. The first and inspiring result in the
multivariate setting was obtained by  E. Kaltofen and P.
Koiran~\cite[Thm~3]{KaKo05}  last year, who showed how to compute
the {\em linear} factors of a   bivariate polynomial $f\in \Q[x,y]$
in  polynomial time in $\ell(f)$. Our result is then an extension of
Kaltofen-Koiran's, and a full generalization of Lenstra's  to the
case of bivariate polynomials.

\bigskip
All  these algorithms (including ours) are based on a {\em gap}
principle first applied by Cucker, Koiran and Smale. The idea is so
strikingly simple and natural that it deserves to be explained. Let
$f\in \Z[x]$ and $\xi\in \Z$ be given, how can we test if
$f(\xi)=0$? Direct evaluation is not feasible, as the size of
$f(\xi)$ can be exponentially big in the input size; an important
exception to this are the easy cases $\xi=0,\pm 1$. Assume that
$f=\sum_{i=1}^ta_ix^{\alpha_i}$ can  be split as
$$f=r+x^u q $$ for  non-zero polynomials $r$ of degree $\deg(r)=k$ and $q$, where
there is a gap between the exponents of $r$ and those of $q$ of
length
$$u-k \ge \log_2 ||f||_1$$
(here $||f||_1:=\sum_{i=1}^t|a_i|$ denotes as usual the {\it
$\ell^1$-norm} of $f$).

 Except for the cases $\xi=0,\pm 1$, this
implies that  $f(\xi)=0$ if and only $q(\xi)=r(\xi)=0$: suppose this
is not the case, namely $f(\xi)=0$ but  $q(\xi)\ne 0$, then
$$
|r(\xi)| \le ||r||_1\cdot |\xi|^{k}  <  ||f||_1\cdot |\xi|^{k} \quad
\mbox{ and } \quad |r(\xi)|= |\xi|^u \cdot |q(\xi)| \ge |\xi|^u$$ so
that $ ||f||_1 >  |\xi|^{u-k} \ge 2^{u-k}$, which contradicts the
gap assumption! Therefore, to test if $f$ vanishes at $\xi \ne 0,\pm
1$, one decomposes $f$ into widely spaced short pieces
$$f=\sum_i x^if_i$$
and tests if $f_i(\xi)=0$ for all $i$.

\medskip
One  crucial fact here is that the decomposition is independent of
the point $\xi$; therefore to find integer roots  it is enough to
find the common roots of a set of low degree polynomials.

The other key ingredient that makes the above argument work is that
any integer $\xi\ne 0,\pm 1$ satisfies a uniform lower bound
$|\xi|\ge2$! In order to apply the same idea to $\xi\in \Q$, the
correct generalization of the absolute value is the {\it height},
defined  as the maximum between numerator and denominator. By
imitating the argument above, but this time for the usual absolute
value {\it and} all the $p$-adic ones, we arrive at the same
conclusion as a consequence that all rational numbers except $0,\pm
1$ have height at least $2$. This is essentially what Lenstra
applied in~\cite{Len99b}; more generally, he was able to handle in
this way other factors besides the linear ones  by considering the
height of their roots after applying a suitable lower bound for
them, namely Dobrowolski's theorem~\cite{Dob79} in the version of
P.~Voutier \cite{Vou96}. In ~\cite{KaKo05}, the authors succeeded to
present  the first  generalization of this gap principle  for
non-univariate polynomials, more precisely for
 linear factors of  bivariate
polynomials.

\medskip
As in these previous works, the key of our algorithm is a suitable
gap theorem. We obtain it as a consequence of a lower bound for the
height of Zariski dense points lying on a curve due to F.~Amoroso
and S.~David~\cite{AmDa00}, as explained in detail in
Section~\ref{The gap theorems}. This result allows to decompose the
given polynomial $f\in \Q[x,y]$ into short pieces; the factors of
$f$ are then computed as the common factors of this low degree
pieces. This strategy works for all factors except the trivial $x$
and $y$ and the  cyclotomic ones, that is, factors which are a
product of binomials (including  monomials) whose coefficients are roots of the unity. As
in the univariate and  linear bivariate cases, these factors have to
be handled separately, see Section~\ref{Computing the low degree
factors of sparse polynomials}.

Since our algorithm operates by reducing to the cases of dense
bivariate and sparse  univariate polynomials, our concern is only to
prove that this reduction can be done in polynomial time in the
sparse encoding. We have not attempted to compute the exponent in
the complexity estimate, which in principle can be quite  big. It is
certainly possible to improve it in view of practical
implementation: in Subsection~\ref{adaptive} we present one idea in
this direction, which consists on adapting the decomposition of $f$
to the size of the candidate factor.

\bigskip As a consequence of  the algorithm, we derive  that the number of
irreducible factors of degree $\le d$ of $f\in \Q[x,y]$  counted
with multiplicities  (different from the trivial factors $x$ or $y$)
 is bounded by $(d\cdot \ell(f))^{O(1)}$.
This is not trivial, as the degree of $f$ can be exponential in
$\ell(f)$, but in fact much better can be said:

\begin{undef}{Proposition 2.}
Let $f\in \Z[x_1,\dots,x_n]$  and consider the factorization
$$
f=q \cdot \prod_p p^{e_p}
$$
where $q$ is a  cyclotomic polynomial, $p\in
\Q[x_1,\dots,x_n]$ runs over all
non-cyclotomic irreducible factors of $f$, and $e_p$
is the corresponding multiplicity, then
$$
\sum_p e_p \le 5^6 \cdot n^3\cdot \log||f||_1 \cdot \log^3(8n\deg(f)).
$$
\end{undef}

In particular the total number of   non-cyclotomic irreducible
factors of any degree of  $f$ is polynomially bounded in terms of
the sparse length of $f$. This fairly unexpected property
generalizes~\cite[Thm~2]{Dob79} and is a further consequence of the
connection with Diophantine Geometry {\it via} the theory of
heights: the Amoroso-David lower bound together with the theorem of
successive algebraic minima of S.-W. Zhang~\cite{Zha95b} imply a
lower bound for the Mahler measure of a non-cyclotomic polynomial,
and from this the statement follows easily.

Moreover, a positive answer to  Lehmer's problem would imply in
the univariate case, see Subsection \ref{plane curves} for details, the stronger estimate
$$
\sum_p e_p \le c \cdot \log||f||_1.
$$
for some absolute constant $c>0$. This is even more surprising,
since it depends on the coefficients of $f$ but not on its degree.
It would be interesting to determine  if it is possible to obtain
such a bound without assuming Lehmer's conjecture.

This should be compared with another result of H.J. Lenstra Jr.: the
total number of irreducible factors of degree $\le d$ of $f\in
\Q[x]$  counted with multiplicities  (different from  $x$)
   is bounded  by
$$ c \cdot t^2\cdot 2^d\cdot  d \cdot \log(2dt)$$
where $t$ is the number of non zero terms of
$f$~\cite[Thm~1]{Len99a}. This bound is exponential, but independent
of the degree and coefficients of $f$. Based on these two results,
it seems natural to consider the following generalization of
Descartes' rule of signs: is the number of all irreducible (and
non-cyclotomic maybe?) factors different from $x$ of a $t$-nomial in
$\Q[x]$  uniformly bounded by some function $B(t)$ depending only on
$t$, and maybe even by $t^{O(1)}$?

\bigskip
Trying to get further, one might ask if it is possible to compute in
sparse polynomial time the {\it absolute} factorization of a sparse
polynomial, that is the irreducible factors over $\Qbarra$. For the
univariate case the answer is clearly ``no'': a univariate
polynomial splits completely as a product of linear factors, and
this cannot be done in sparse polynomial time. For the bivariate
case, it can be shown that the computation of binomial factors  is
equivalent to the factorization of a univariate polynomial, so that
binomials factors over $\Qbarra$ cannot be computed either.

Here, we show that except for these, we can compute all other
irreducible factors over $\Qbarra$ of low degree, in sparse
polynomial time.  To give sense to such a statement, we have to
specify the way algebraic coefficients are handled: a number field
$K$ is described by an irreducible monic polynomial
$g=\sum_{j=0}^{\delta-1} g_jz^j\in \Z[z]$ such that $K=\Q(\theta)$
for one of its roots, and this $g$ is given in dense representation
 by the list of all coefficients  $g_j$ in some specified order, including the zero ones.
Each irreducible factor $p$ in the output of the algorithm is
encoded by giving a number field $K$ such that $p\in K[x,y]$ and by
the dense list of its coefficients, each coefficient $b\in K$ being
represented by its vector of rational components
$b:=(b_{0},\dots,b_{\delta-1})$ with respect to the basis
$(\theta^j)_{0 \le j\le  \delta-1}$.

\begin{undef}{Theorem 3.}
There is a deterministic algorithm that, given $f\in \Q[x,y]$ and $d
\ge 1$, computes all irreducible factors of $f$ in $\Qbarra[x,y]$ of
degree $\le d$, together with their multiplicities, except for the
binomial ones, in $(d\cdot\ell(f))^{O(1)}$ bit operations.
\end{undef}

This algorithm follows from another suitable gap theorem that we obtain as
a consequence of a further result of Amoroso and David, a quantitative version of
the Bomogolov problem over the torus~\cite{AmDa03}.
Furthermore, we deduce from
their result an estimate for the number of
non-binomial factors of a given
$f\in \Qbarra[x_1,\dots, x_n]$ (Proposition~\ref{factores2}).

\bigskip Several interesting questions arose during our  work.
The most obvious is the extension of these algorithms to
multivariate polynomials; this seems quite feasible as the necessary
lower bounds for the height of points in a hypersurface  already
appeared in the literature~\cite{AmDa00,AmDa03,Pon01,Pon05b}.

An interesting open problem is the following: the restriction to
computing bounded degree factors keeps their length under control,
giving the possibility of computing them in sparse polynomial time.
But, what if we look for factors with a fixed number of monomials,
can we still find all of them in sparse polynomial time? For
instance, can we compute all trinomial factors
$$
p=a_1x^{\alpha_1}+a_2x^{\alpha_2}+a_3x^{\alpha_3} \in \Q[x]
$$
of a given $f\in \Q[x]$ in polynomial time?

\bigskip The outline of the paper is as follows.
In Section~\ref{Heights} we explain the basics of the height theory for
points, polynomials and curves, and we prove the upper bounds for the number
 of factors of a sparse polynomial.
In Section~\ref{The gap theorems} we obtain the gap theorems,
as a consequence of the lower bounds for the height of points on curves.
In Section~\ref{Computing the low degree factors
of sparse polynomials} we present the algorithms for rational and absolute factorization
 and estimate their theoretical complexity.

\bigskip
{\bf Acknowledgements}

We thank Corentin Pontreau for helpful discussions on lower bounds
for the height.

The core of this paper was written during October--December 2005
while M.~Sombra was visiting the University of Buenos Aires,
Argentina; he particularly thanks Ricardo Dur\'an for his
invitation. He also thanks the Mathematical Sciences Research
Institute at Berkeley, USA, where he stayed during January 2006.

\bigskip
{\bf Note}

We learned on January 23th, 2006,  that a multivariate version of
Theorem~1 had independently been achieved by Erich Kaltofen and
Pascal Koiran. We immediately sent them the present paper, which  at
that time was in essentially final form.

\section{Heights} \label{Heights} \label{ Height of polynomials over }

Throughout this paper $\Q$ denotes the field of rational numbers,
$K$ a number field,  $L$ a finite extension of $K$, $\Qbarra$ an
algebraic closure of $\Q$ and  $G_\infty$ the subset of $\Qbarra$ of all roots of
the unity.  We denote by $\A^n$ the affine space of
$n$ dimensions over $\Qbarra$. For a polynomial $p\in
\Qbarra[x_1,\dots,x_n]$ we denote by $Z(p)\subset \A^n$ the affine
hypersurface defined by $p$. A curve or a variety is assumed to be
equidimensional; by irreducibility of a variety we understand its
geometric irreducibility, that is with respect to $\Qbarra$.

\smallskip
For every rational prime $p$  we denote by $ | \cdot |_p $ the
$p$-adic absolute value over $\Q$ such that
$|p|_p=p^{-1}$.  We also denote the ordinary absolute value over
$\Q$ by $|\cdot|_\infty$ or simply by $ | \cdot | $. These form a
complete set of independent absolute values over $\Q$: we identify
the set $M_\Q$ of these absolute values with the set $\{ \infty, p
\,;\, \, p \ \mbox{prime} \}  $.
More generally, we write $M_K$ for the set of absolute values over $K$ extending
the absolute values in $M_\Q$,
and we note by $M_K^\infty$ the subset
of Archimedean absolute values of $M_K$.

\smallskip
For $v_0\in M_K $ we denote by $\Q_{v_0}$ the completion of $\Q$ with
respect to the absolute value $v_0$.
In case $v_0=\infty$ we have
$\Q_\infty=\R$, while in case $v_0=p$ is a prime, $\Q_p$ is
the $p$-adic field.
There exists a unique extension of $v_0$ to an
absolute value over the algebraic closure $\overline{\Q}_v$.
For $v\in M_K $ we also denote by $K_v$ the completion of $K$ with
respect to $v$.
If $v$ extends an absolute value $v_0\in M_\Q$, then
$K_v$ is a finite extension of $\Q_{v_0}$.
We denote $\sigma_v:K\hookrightarrow \Qbarra_v$
a (not necessarily unique) embedding
 corresponding to $v$, that is such that
$|a |_v=|\sigma_v(a)|_{v_0}$ for every $a \in K$.

\subsection{Height of points and polynomials} \label{algheight}

In this subsection we introduce the basic definitions and properties of the
height of points and polynomials that we will use in the sequel.
We refer for instance to
\cite{HS00} for a complete treatment.

\medskip
 The  {\it
(logarithmic) height} $h(\xi)$ of an algebraic number $\xi\in \Qbarra$ can be defined in terms of its
primitive integer minimal polynomial
$$
p_\xi(x)= c\cdot
 \prod_{\sigma: K\hookrightarrow \Qbarra}
(x-\sigma(\xi))  \quad \in \Z[x]
$$
where $\sigma$ runs over all $\Q$-embeddings of $K:=\Q(\xi)$ in $\Qbarra$,
by the formula
\begin{equation}\label{mahler}
h(\xi)=\frac{1}{[K:\Q]} \left(
\log |c| +\sum_{\sigma: K\hookrightarrow \Qbarra} \max\{0,\log|\sigma(\xi)|\}\right).
\end{equation}
We have $h(\xi)\ge 0$, and $h(\xi)=0$ if and only if either $\xi=0$
or $ \xi\in G_\infty$, the subset of $\Qbarra$ of all roots of~1 (Kronecker's theorem).
Besides, for a rational  $\xi=m/n\in \Q^\times$ in reduced expression, we
easily check that
$h(\xi)=\log\max\{|m|,n\}$.
Alternatively,
the height can be defined {\it via} the Mahler measure
of the minimal polynomial as
$$
m(p_\xi):=\int_0^1 \log|p_\xi(e^{2\pi i u})|\,du= [K:\Q] \cdot h(\xi);
$$
this identity is a consequence of Jensen's formula.

\medskip
More generally, the {\it height} of a point
$\xi:=(\xi_1,\dots,\xi_n)\in \A^n$ is defined {\it via} the Weil
formula
$$
h(\xi) := {1\over [K:\Q]}\sum_{v\in M_K} {[K_v:\Q_v]}
\log \max \{ 1,|\xi_1|_v,\dots,|\xi_n|_v\}
$$
for any number field $K$ containing the coordinates $\xi_i$.
For $n=1$ this gives
$$
h(\xi) = {1\over [K:\Q]}\sum_{v\in M_K} {[K_v:\Q_v]}
\log \max \{ 1,|\xi|_v\}
$$
and it can be shown that this coincides with the previous definition.
With this expression we readily verify that
for $\xi, \eta\in \Qbarra $ we have that  $h(\xi\cdot \eta)\le h(\xi)+h(\eta)$ and
$$
h(\xi^n)=|n|\,h(\xi) \quad \mbox{\rm  for } n\in\Z;
$$
in particular
$h(\xi^{-1})=h(\xi)$ and $h(\omega\cdot \xi) =h(\xi)$ for any
root of unity $\omega\in G_\infty$.
We will be mostly interested on points in the plane
$\xi=(\xi_1,\xi_2)\in \A^2$,
in that case the formula reduces to
$$
h(\xi) = {1\over [K:\Q]}\sum_{v\in M_K} {[K_v:\Q_v]}
\log \max \{ 1,|\xi_1|_v,|\xi_2|_v\}.
$$

\bigskip
Now we introduce a few notions for the height of a polynomial
that will prove useful in the sequel.
We will restrict to bivariate polynomials, although it is clear
that all this extends to the multivariate case.

\medskip
For a polynomial
$f = \sum_{i=1}^t a_i \, x^{\alpha_i}y^{\beta_i} \in K[x,y]$,
its
{\it absolute value} with respect to $v\in M_K$ is
$$
|f|_v:=\max \{|a_1|_v,\dots,|a_t|_v\}.
$$
The {\it height} of $f$ is then defined as
$$
h(f) := {1\over [K:\Q]}\sum_{v\in M_K} {[K_v:\Q_v]}
\log (|f|_v),
$$
which is invariant by scalar multiplication because of the product formula
$$\sum_{v\in M_K} [K_v:\Q_v]
\log (|a|_v)=0 ,\ \ \forall \ a\in K^\times.$$

Therefore $h(f)$ is the Weil height of the projective point $(a_1:\cdots:a_t)$.
This is independent of the chosen field $K$ as long as it contains
all of the $a_i$'s.

\medskip
For a bivariate polynomial with complex coefficients $f\in \C[x,y]$
we consider the {\it Mahler measure}
$$
m(f):=\int_0^1\int_0^1 \log|f(e^{2\pi i u},e^{2\pi i v})|\,d
u\,d v,
$$
and for a polynomial $f\in K[x,y]$ with {\it algebraic} coefficients we define its
{\it (global)} {\it Mahler measure} by the adelic formula
$$
m_\Qbarra(f):= {1\over [K:\Q]}\left(
 \sum_{v\in M_K^\infty} [K_v:\Q_v] \, m(\sigma_v(f)) +
\quad \sum_{v\notin M_K^\infty}[K_v:\Q_v]\log|f|_v\right).
$$

\medskip
We also consider the height associated to the $\ell^1$-norm:
$$
h_1(f):={1\over [K:\Q]}\left(
\sum_{v\in M_K^\infty} [K_v:\Q_v]\, \log||\sigma_v(f)||_1 +
\quad \sum_{v\notin M_K^\infty} [K_v:\Q_v]\log|f|_v\right).
$$

\medskip
For a {\it primitive} $f\in \Z[x,y]$, these notions give
$$
h(f) = \log |f| =\log \max \{|a_1|,\dots,|a_t|\},
\
h_1(f) =
\log ||f||_1 =\log (|a_1|+ \cdots + |a_t|), \
m_\Qbarra(f)=m(f).
$$
All these are
invariant by scalar multiplication.
In general for any $f\in \Q[x,y]$
write $f=c \cdot \widetilde f$ for some
$c\in \Q^\times$ and
$\widetilde f\in \Z[x,y]$ the
 primitive polynomial
with integer coefficients associated to $f$, then
$h(f)= \log |\widetilde f|$, $h_1(f)=\log||\widetilde f||_1$ and
$m_\Qbarra(f)=m(\widetilde f)$.

\medskip
We will use the following comparison between the heights of
a given $f \in K[x,y]$, which can be directly proven from the definitions:
\begin{equation} \label{comparaciones}
h(f), m_\Qbarra(f) \le h_1(f)\le h(f) + \log(t).
\end{equation}

\subsection{Height {\it of} and {\it on} plane curves} \label{plane curves}

A plane curve $C\subset \A^2$ can have some isolated  points of small height.
For instance the line
$$
Z(x+y-1) \subset \A^2
$$
has the points $(1,0), (0,1), \big((1\pm\sqrt{3})/2, (1\mp\sqrt{3})/2\big)$ all of
whose coordinates are roots of 1 and so their height is 0.
D. Zagier~\cite{Zag93} showed that the height of any other point $\xi\in Z(x+y-1)$
is bounded from below by a positive constant
$$
h(\xi) \ge h(\xi_0) = 0.1911
$$
where $\xi_0$ denotes the largest real root of the polynomial $x^6-x^4-1$.
Somehow the fact that a curve has some torsion points on it does not
reflect  its general behavior.
A more interesting parameter is the height of a Zariski dense set of points.
This is measured by the {\it essential minimum}, which
for a plane curve $C\subset \A^2$ is defined as
$$
\mu^\ess(C) :=
\inf \big\{\eta \ge 0 :  \{ \xi \in C: h(\xi) \le \eta\}
\mbox{ is an infinite set } \big\}.
$$
For instance, thanks to Zagier's result,
$$
\mu^\ess(Z(x+y-1))\ge 0.1911.
$$
This is a particular case of
the Bogomolov problem over the torus proved by
Zhang~\cite{Zha95b}
which  asserts that for a subvariety of $\T^n:=(\Qbarra^\times)^n$,
the vanishing of the essential minimum is equivalent to being
{\it torsion}. This result, and others we are going to use, are stated for
the torus, but $\T^n$ is naturally embedded as an open subset of $\A^n$,
and since these results depend on Zariski dense sets,  they can all be translated to $\A^n$.

For an irreducible plane curve $C\subset \A^2$, being {\it torsion} is
equivalent to say that  there exist
$\alpha,\beta\ge 0$ not both zero, and  $\omega\in G_\infty\cup \{0\}$
such that
$$
\mbox{ either } \quad  C= Z(x^\alpha - \omega y^\beta)
\quad \mbox{ or } \quad
C= Z(x^\alpha y^\beta-\omega).
$$
The irreducible curve $C$ is (we should rather say ``corresponds to'') a
{\it translate of a subgroup} whenever
there exists  $\xi\in \Qbarra$ such that
$$
\mbox{ either } \quad  C= Z(x^\alpha - \xi y^\beta)
\quad \mbox{ or } \quad
C= Z(x^\alpha y^\beta-\xi).
$$
By definition, a general affine plane curve is torsion
(resp. translate of a subgroup) if and only if all  its irreducible
components are so. The statement of the Bogomolov
problem (now a theorem)
is that $\mu^\ess(C)=0$ if and only if $C$ is torsion. In other words,
if $C$ is not of this form, there exists a positive constant
$c(C)>0$  such that
$$
h(\xi) \ge c(C) \quad \mbox{ for all but a finite number of } \xi \in C.
$$

\medskip
There is an extension of the notion of Weil height of points to higher-dimensional varieties.
This notion was first introduced by P. Philippon \cite{Phi91};
for an irreducible hypersurface $V\subset\A^n$ defined by a
polynomial $p\in K[x_1,\dots,x_n]$,
it coincides with the global Mahler measure of $p$ \cite{DaPh99,Pon01}:
\begin{equation}\label{h y m}
h(V)= m_\Qbarra(p).
\end{equation}

\bigskip
The distribution of the height of algebraic points in a curve is in
close connection with the height of the curve itself. The relation
is given by the {\it theorem of algebraic successive minima} of
Zhang~\cite[Thm~5.2 and Lem.~6.5(3)]{Zha95b}:
$$
\mu^{\ess}(C)\le \frac{h(C)}{\deg(C)}\le 2  \mu^{\ess}(C).
$$
Actually, Zhang's result is more precise (all successive
minima appear, not only the first one which is  the essential minimum) and more
general, as it works for varieties of any dimension and
for any ``reasonable'' height function.

The stated version is sufficient for our application; for a more elementary
proof we refer to~\cite[\S~6]{DaPh99}.
It is an open problem to determine if this estimate is optimal for the case of plane
 curves or more generally for hypersurfaces (
it has been shown to be optimal if we allow varieties of higher
codimension~\cite[Thm~5.1]{PS04}). Thanks to this result, the
Bomogolov problem for plane curves can be rephrased as  $h(C)=0$ if
and only if $C$ is torsion. Under this form, the conjecture was
already proven by W.~Lawton in 1977~\cite{Law77}.

\bigskip
For $\xi \in \Qbarra^\times$ we have that $h(\xi)=0$ if and only if
$\xi\in G_\infty$; this is the 0-dimensional (easy) case of the Bogomolov
problem.
Lehmer's conjecture gives a lower bound for the height of non-torsion points,
its statement being that there exists a positive constant $c>0$
such that
$$
h( \xi) \ge \frac{c}{[\Q(\xi):\Q]} \quad \mbox{\rm for } \xi \notin G_\infty.
$$
This conjecture has been widely generalized.
Here we are only interested in the case of curves:

\begin{conjecture} \label{lehmer+bogomolov}
\

\begin{itemize}
\item[(i)]
{\em Lehmer problem for plane curves:}
let $C\subset \A^2$ be an irreducible curve defined over a number field $K$
which is not torsion,
then there exists a universal $c>0$ such that
$$
\mu^\ess(C) \ge \frac{c}{[K:\Q]\deg(C)}.
$$
\item[(ii)] {\em Effective Bogomolov problem for plane curves:}
let $C\subset \A^2$ be an irreducible curve which is not a translate
of a subgroup,
then there exists a universal $c>0$ such that
$$
\mu^\ess(C) \ge \frac{c}{\deg(C)}.
$$
\end{itemize}
\end{conjecture}

These two conjecture look similar but they are not. The
generalization of  Lehmer problem is of arithmetic nature since the
degree of the number field plays a role, while the quantitative
Bogomolov problem is of  geometric nature since it makes no
reference to the field of definition. It has been shown that
conjecture~\ref{lehmer+bogomolov}(i) is implied by the classical
Lehmer problem~\cite{Law77}. Conjecture~\ref{lehmer+bogomolov}(ii)
is~\cite[Conj.~1.1]{DaPh99}.

Because of the theorem of successive minima, it is equivalent to
have lower bounds for the essential minimum or for the height,
that is the (global) Mahler
measure of the defining polynomial of $C$.

\medskip
Nowadays all these results are proved ``up to an $\varepsilon$'':
for the Lehmer problem we will be mainly applying the following lower
bound due to
Amoroso and David~\cite{AmDa00},  in the version of
C. Pontreau~\cite[Prop.~IV.1]{Pon05}
who simplified the proof and made all constants explicit:
let  $C \subset \A^2$ be a non-torsion curve defined by an irreducible
polynomial $p\in
\Z[x,y]$ and set $d:=\deg(C)=\deg(p)$, then
  \begin{equation}\label{cotaPontreau}
  \mu^\ess (C)\ge \frac{1}{5^6d}\times \left(\frac{\log \log (16d)}{\log
  (16d)}\right)^3.
\end{equation}
In the reference this result is stated in terms of $h(C)$; you have
to look into the proof for the version up here. In fact we will be
using the version over a number field:

\begin{cor}\label{lo que usamos}
Let $C\in \A^2$ be a curve
defined by an irreducible polynomial $p\in
K[x,y]$ which is not
of the form $ p=\prod_i (x^{\alpha}
 -\omega_i y^{\beta})$ nor $ p=\prod_i(x^{\alpha}y^{\beta}
 -\omega_i) $ for some $\alpha,\beta\ge 0$ not both zero and
$\omega_i\in G_\infty \cup \{0\}$ and set $d:= \deg(C)=\deg(p)$, then
$$
\mu^\ess (C)\ge  \frac{1}{5^6[K:\Q] d}
\times \left(\frac{\log \log (16[K:\Q] d)}{\log (16[K:\Q] d)}
\right)^3.$$
\end{cor}

This follows immediately from~(\ref{cotaPontreau}) by considering
the norm $N(p):=\displaystyle{\prod_{\sigma:K\hookrightarrow
\Qbarra} \sigma(p)} \in \Q[x,y]$.

\medskip
For the effective Bogomolov problem we use another result of Amoroso
and David: let $C \subset \A^2$ be a curve which is not a translate
of a subgroup and $d:=\deg(C)=\deg(p)$, then~\cite[Thm~1.5]{AmDa03}:
\begin{equation} \label{effective_bogomolov}
\mu^\ess (C)\ge \frac{1}{2^{70} d}\times \frac{(\log \log (d+2))^4}{(\log
(d+2))^5}.
\end{equation}

\subsection{On the number of factors of a sparse polynomial}

General lower bounds for the Mahler measure yield immediately upper
bounds for the number of factors of a given polynomial. To our
knowledge, this observation appears for the first time in the work
of E.~Dobrowolski~\cite{Dob79}. Here we treat the general
$n$-dimensional case. The notions and results of the previous
subsection extend to hypersurfaces. We will state them but instead
refer the interested reader to the literature for $n\ge 3$.

We recall that a polynomial is cyclotomic if it is a product of binomials (including monomials)
whose coefficients are roots of the unity.

\begin{prop} \label{factores1}
Let  $f\in K[x_1,\dots,x_n]$ and consider the factorization
$$
f=q \cdot \prod_p p^{e_p}
$$
where $q$ is cyclotomic, $p\in K[x_1,\dots,x_n]$  runs over
all non-cyclotomic
irreducible factors of $f$, and $e_p$  is the corresponding
multiplicity, then
$$
\sum_p e_p \le 5^6 \cdot n^3 \cdot [K:\Q] \cdot h_1(f) \cdot
\log^3(8n[K:\Q]\deg(f)).
$$
\end{prop}

\begin{proof}
We have that $m_\Qbarra(q)=0$ as $q$ is cyclotomic an so
$$
\sum_p e_p m_\Qbarra(p) = m_\Qbarra(f) \le h_1(f).
$$
For each non-cyclotomic factor $p\in \Q[x_1,\dots,x_n]$ we minorate
the Mahler measure by the Amoroso-David's lower bound in the version
of Pontreau~\cite[Thm~1.6]{Pon01} (see the
estimate~(\ref{cotaPontreau}) above for the case $n=2$), from which
we derive that if $V\subset \A^n$ is an hypersurface defined by an
irreducible polynomial over $K$, then
$$[K:\Q]\cdot h(V)\ge \frac{1}{5^6 \cdot n^3}\cdot \left(\frac{\log(n\log(8n[K:\Q]
\deg(V))}{\log(8n[K:\Q] \deg(V))}\right)^3.$$

Therefore, by Identity (\ref{h y m}), we have
$$
 [K:\Q]\cdot m_\Qbarra(p)
\ge
\frac{1}{5^6 \cdot n^3 \cdot  \log^3(8n[K:\Q] \deg(p))}  \ge
\frac{1}{5^6 \cdot n^3 \cdot  \log^3(8n[K:\Q] \deg(f))},
$$
which implies
$$
[K:\Q] \cdot h_1(f) \ge \frac{1}{5^6 \cdot n^3
\cdot  \log^3(8n[K:\Q] \deg(f))} \, \sum_p e_p
$$
from where we deduce our result.
\end{proof}

This is a generalization to $n\ge 2$ of~\cite[Thm~2]{Dob79}. As
said, a positive answer to the classical Lehmer problem would imply
a positive lower bound for the Mahler measure of an arbitrary
non-cyclotomic polynomial $p\in K[x_1,\dots,x_n]$, of the form
$$
m_\Qbarra(p)\ge \frac{c}{[K:\Q]}
$$
for some universal constant $c>0$, namely
Conjecture~\ref{lehmer+bogomolov}(i).
Applying this to the argument above, the previous proposition would improve to
\begin{equation} \label{lehmer}
\sum_p e_p \le c^{-1} \cdot [K:\Q] \cdot h_1(f).
\end{equation}
In a similar way, we can produce an upper bound for the number
of non-binomial irreducible factors over $\Qbarra$:

\begin{prop} \label{factores2}
Let $f\in \Qbarra[x_1,\dots,x_n]$ and consider the factorization
$$
f=q \cdot \prod_p p^{e_p}
$$
were $q$ is a product of binomials, $p\in \Qbarra[x_1,\dots,x_n]$  runs over
all non-binomial
irreducible factors of $f$, and $e_p$  is the corresponding
multiplicity, then
$$
\sum_p e_p \le 10^{14} \cdot n^8 \cdot h_1(f) \cdot
\log^5 (\max \{16, n\deg(f)\}).
$$
\end{prop}

\begin{proof}
We have that
$$
\sum_p e_p m_\Qbarra(p) = m_\Qbarra(f) \le h_1(f);
$$
apply the Amoroso-David quantitative Bogomolov problem in the
version of Pontreau~\cite[Thm~1.5]{Pon05b}
(or~(\ref{effective_bogomolov}) above for the case $n=2$).
\end{proof}

Similarly, a positive answer to the effective Bogomolov problem
(Conjecture~\ref{lehmer+bogomolov}(ii))
would imply that
$$
\sum_p e_p \le c^{-1} \cdot h_1(f)
\quad  \mbox{ for a universal constant } c>0.
$$

\section{Gap theorems} \label{The gap theorems}

By a {\it gap theorem}, following \cite{CKS99,Len99b,KaKo05}, we understand a statement asserting
that for a polynomial $f$ decomposed as
$$f=r+s$$
for non-zero polynomials $r$ and $s$, then $f$ has a given property
if and only $r$ and $s$ have it, provided that $r$ and $s$ are
sufficiently separated. We introduce some notation:

\begin{defn} For $p\in \overline\Q[x,y]$ such that
$\deg_y(p)\ge 1$ we set
$$\lambda(p):=\inf \big\{\eta \ge 0 :  \{ (\omega,\nu) \in G_\infty \times \Qbarra
:  p(\omega,\nu)=0, \ h(\nu)\le \eta\} \mbox{ is an infinite set }\big\}.
$$
\end{defn}

Since $\deg_y(p)\ge 1$, for all but a finite number of
 $\omega\in G_\infty$
there exists some $\nu\in \Qbarra$ such that
$p(\omega,\nu)=0$ and so $\lambda(p)$ is well-defined and non-negative.

In what follows we deal with irreducible polynomials, that
 are defined up to a scalar factor. For simplicity
we always refer to  one (obvious) representant in each class of
associate irreducible polynomials.

The following is the main result of this section:

\begin{teo} \label{gapMartin}
Let $f,r,q \in \Qbarra[x,y]$ be such that $ f=r+y^{u}\cdot q$. Let
also be given an irreducible polynomial $p\in \Qbarra[x,y]$,
 $p\ne y$, such that $\deg_y(p)\ge 1$, and suppose that
$$
(u-\deg_y(r))\cdot \lambda(p) \ge
h_1(f),$$
then $p$ divides $f$ if and only if it divides $r$ and $q$.
\end{teo}

For its proof we need the following lemma:

\begin{lem} \label{prep}
Let $f,r,q \in \Qbarra[x,y]$ be such that
$
f=r+y^{u}\cdot q$.
Let also be given
$\omega\in
G_\infty $  and $\nu\in \overline \Q^\times$ be such that $f(\omega,\nu)=0$
but
 $q(\omega, \nu)\ne 0$, then there exists a constant $\delta(f) >0$ not
depending on
$(\omega,\nu)$ such that
$$(u-\deg_y(r))\cdot h(\nu) \le h_1(f)-\delta(f).$$
\end{lem}

\begin{proof}
Let  $K$ be  a number field containing
the coefficients of $f$,
$\omega$ and $\nu$, and set $k:=\deg_y(r)$.
For each absolute value $v\in M_K$ we have two cases:
\begin{itemize}
\item $|\nu|_v\le 1$:  since $|\omega|_v=1$ we have that
$$
|q(\omega,\nu)|_v\le \left\{ \begin{array}{ll}
||\sigma_v(q)||_{1} & \quad \mbox{\rm for } v\in M_K^\infty ,\\[2mm]
|q|_{v} & \quad \mbox{\rm for } v\notin M_K^\infty .
\end{array} \right.
$$
\item $|\nu|_v>1$:  using that
$f(\omega,\nu)=r(\omega,\nu)+\nu^{u}q(\omega,\nu)=0$ we infer that
$$
|\nu|_v^{u}\cdot |q(\omega,\nu)|_v=|r(\omega,\nu)|_v
\le \left\{ \begin{array}{ll}
 |\nu|_v^{k} \cdot ||\sigma_v(r)||_{1} & \quad \mbox{\rm for } v\in M_K^\infty ,\\[2mm]
 |\nu|_v^{k} \cdot |r|_{v}  & \quad \mbox{\rm for } v\notin M_K^\infty .
\end{array} \right.
$$
\end{itemize}
As both $r$ and $q$ are non-zero,
$||\sigma_v(q)||_{1}, ||\sigma_v(r)||_{1} < ||\sigma_v(f)||_{1}$
and so
$$
\log ||\sigma_v(q)||_{1},  \ \log ||\sigma_v(r)||_{1} \le
\log ||\sigma_v(f)||_{1} - \delta(f)
$$
for some $\delta(f) >0$ depending only on
$f$.
The previous inequalities imply that
$$
(u-k) \log\max\{1,|\nu|_v\} +
  \log|q(\omega,\nu)|_v
\le \left\{ \begin{array}{ll}
 \log ||\sigma_v(f)||_{1} - \delta(f)
& \quad \mbox{\rm for } v\in M_K^\infty ,\\[2mm]
 \log|f|_v& \quad \mbox{\rm for } v\notin M_K^\infty .
\end{array} \right.
$$
By summing up over all absolute values, using the product formula
and the definition of the height, one obtains that
\begin{align*}
(u-k)\cdot h(\nu) &={1\over [K:\Q]}  \sum_{v\in M_K}
[K_v:\Q_v]\,\big((u-k) \log\max\{1,|\nu|_v\} +
  \log|q(\omega,\nu)|_v \big) \\[2mm]
&\le {1\over [K:\Q]} \left(\sum_{v\in M_K^\infty} [K_v:\Q_v] \big(
\log ||\sigma_v(f)||_{1} - \delta(f)\big) \quad +  \sum_{v\notin
M_K^\infty} [K_v:\Q_v]
 \log |f|_{v} \right)\\[2mm]
&=h_1(f)-\delta(f).
\end{align*}
\end{proof}

\begin{Proof}{Proof  of Theorem~\ref{gapMartin}}
The ``$\Leftarrow$'' is trivial, so we show the other implication.

Suppose that $p\div f$ but $p\ndiv q$. From the fact that $p$ is
irreducible we have that the set of common roots of $p$ and $q$ is
finite. Also, since $\deg_y(p)\ge 1$ and $p\ne y$, the set
$\{(\omega,\nu)\in G_\infty\times \Qbarra^\times: p(\omega,\nu)=0\}$
is infinite. Given $\varepsilon >0$, it follows from the definition
of $\lambda(p)$ that the set $\{(\omega,\nu)\in G_\infty\times
\Qbarra: p(\omega,\nu)=0,\ h(\nu)\le \lambda(p)-\varepsilon\}$ is
finite. Therefore there exist an infinite number of $(\omega,\nu)\in
G_\infty\times \Qbarra^\times$ such that $p(\omega,\nu)=0$ and
$h(\nu)> \lambda(p)-\varepsilon$, and there still exist
 some $\omega\in G_\infty$ and
$\nu \in \Qbarra^\times$ such that
$$
p(\omega,\nu)=0, \quad q(\omega,\nu)\ne 0  \quad \mbox{\rm and }
h(\nu)> \lambda(p)-\varepsilon.
$$
Applying   Lemma~\ref{prep}
$$ (u-k) \, (\lambda(p)-\varepsilon)\le (u-k) \, h(\nu)\le h_1(f)-\delta(f).$$
Since this holds for all $\varepsilon>0$, we infer
$$
(u-k)\,\lambda(p)\le h_1(f)-\delta(f) <h_1(f)
$$
because $\delta(f)$ does not depend on $(\omega,\nu)$ and so does not depend on $\varepsilon$ either.
This contradicts the hypothesis: $(u-k)\,\lambda(p)\ge h_1(f)$.
Therefore $p\div q$ and $p\div -y^u\cdot q= r$ as wanted.
\end{Proof}

\bigskip
Of course this result is only useful whenever
$\lambda(p)>0$. What happens is that this parameter is bounded from below
by the essential minimum, and so all existing estimations for the essential minimum
will give us a corresponding gap theorem.

\begin{lem} \label{comparacion}
Let  $p$ be an irreducible polynomial in $K[x,y]$ such that
$\deg_y(p) \ge 1$. Then $$ \lambda(p) \ge \mu^\ess(Z(p)).$$
\end{lem}

\begin{proof}
Observe that $h(\nu)=h(\omega,\nu)$; we can then rephrase the definition
of $\lambda(p)$ as
$$
\lambda(p) =
\inf \big\{\eta \ge 0 :  \{ \xi \in Z(p) \cap (G_\infty \times \Qbarra)
: h(\xi) \le \eta\}
\mbox{ is an infinite set } \big\}.
$$
Compare with the definition of the essential minimum:
$$
\mu^\ess(Z(p)) =
\inf \big\{\eta \ge 0 :  \{ \xi \in Z(p): h(\xi) \le \eta\}
\mbox{ is an infinite set } \big\},
$$
so that $\lambda(p)$ is the infimum over a subset of the set used to
define $\mu^\ess(Z(p))$  and the inequality is clear.
\end{proof}

Equality in Lemma~\ref{comparacion} above does not necessarily hold:
consider $p:=x^\alpha-\xi y^\beta$, then for any $(\omega,\nu)\in
G_\infty\times \Qbarra$ we have that  $p(\omega,\nu)=0\iff \nu^\beta=\omega^\alpha/ \xi $
and so
$$
h(\nu)=\frac{h(\nu^\beta)}{\beta}= \frac{h(\omega^\alpha/ \xi )}{\beta}=
\frac{h(\xi)}{\beta}.$$
Hence
$$
\lambda(p(x,y))=h(\xi)/\beta \quad \mbox{ while } \quad
\lambda(p(y,x))=h(\xi)/\alpha,
$$
in particular  $\lambda$
depends on the order of the variables, while of course the essential
minimum does not, so there cannot coincide in general.
One can prove that $\mu^\ess(p)=h(\xi)/\max\{\alpha,\beta\}$~\cite[Prop.~5.4]{PS04}.

\bigskip
From Corollary~\ref{lo que usamos} we deduce:

\begin{cor} \label{gap_aritmetico}
Let $f,r,q \in K[x,y]$ be such that $ f=r+y^{u}\cdot q$. Let also be
given an irreducible $p\in K[x,y]$
 that is not
of the form $ p=\prod_i (x^{\alpha}
 -\omega_i y^{\beta})$ nor $ p=\prod_i(x^{\alpha}y^{\beta}
 -\omega_i) $ for some $\alpha,\beta\ge 0$ not both zero and
$\omega_i\in G_\infty \cup \{0\}$, and set $d:= \deg(p)$. Suppose that
$$
u-\deg_y(r)\ge
5^6\cdot[K:\Q] \cdot d \cdot \log^3 (16[K:\Q] d)\cdot
h_1(f),
$$
then $p$ divides $f$ if and only if it divides $r$ and $q$.
\end{cor}

Similarly we obtain the following gap theorem from the lower bound~(\ref{effective_bogomolov}):

\begin{cor} \label{gap_geometrico}
Let $f,r,q \in \Qbarra[x,y]$ be such that $ f=r+y^{u}\cdot q$. Let
also be given an irreducible $p\in \Qbarra[x,y]$ which is not a
binomial, and set $d:= \deg(p)$. Suppose that
$$
u-\deg_y(r)\ge
2^{70} \cdot d \cdot \log^5(d+2) \cdot h_1(f),
$$
then $p$ divides $f$ if and only if it divides $r$ and $q$.
\end{cor}

\section{Computing the low degree factors of sparse polynomials}
\label{Computing the low degree factors of sparse polynomials}

The goal of this section is to present the rational and absolute
factorization algorithms for sparse bivariate polynomials. Our
conventions about encoding are the usual ones, the same as in for
instance~\cite{Len99b}. The number of bits needed to write down an
integer $a\in \Z$ is $\lfloor \log_2 (a) \rfloor +1$ for the digits
and 1 more for the sign. For a rational $a=m/n \in \Q$ in reduced
expression, we define its bit length as
$$
\ell(a)= \ell(m)+\ell(n)-2 = \lfloor \log_2 |m|\rfloor+\lfloor
\log_2 (n) \rfloor+2;
$$
the somewhat artificial ``$-2$'' is there just to make this coincide
with the previous notation for an integer $a$. The sparse encoding
of $f=\sum_{i=1}^t a_ix^{\alpha_i}y^{\beta_i}\in \Q[x,y]$ is the
list $(a_i,\alpha_i,\beta_i)_{1\le i\le t}$ of its (non-zero)
coefficients and corresponding exponents, and so its {bit length} is
\begin{equation} \label{ell(f)}
\ell(f):=\sum_{i=1}^t \Big(\ell(a_i) + \lfloor \log_2 (\alpha_i) \rfloor
+\lfloor \log_2(\beta_i) \rfloor+2\Big);
\end{equation}
observe that $\ell(f)$ is an upper bound for
$ t $, $\log_2(\deg f)$ and $h(f)$, and in fact is polynomially
equivalent to these quantities:
$\ell(f)=(t\cdot \log_2(\deg f) \cdot h(f))^{O(1)}$.

\medskip
For encoding polynomials over number fields we have to say how
number fields and algebraic numbers are handled: a number field $K$
of degree $\delta=[K:\Q]$ is described  by an irreducible monic
polynomial $g=\sum_{j=0}^{\delta-1} g_jz^j\in \Z[z]$ such that
$K=\Q(\theta)$ for one of its roots, and this $g$  is given in dense
representation by the (ordered) list of all its coefficients $g_j$
including the zero ones. The length of this description is
$$
\ell(K) := \sum_{j=0}^{\delta-1} \ell(g_j);
$$
in particular $ \ell(K) \ge [K:\Q], h(g)$.
An element  $b\in K$ is represented by
its vector of rational components $(b_{0},\dots,b_{\delta-1})$ with respect
to the basis $(\theta^j)_{0\le j\le \delta-1}$.
It can be shown by (you need some estimate between the height of an algebraic integer and that
 of its minimal polynomial) that
$$
h(b) \le \ell_K(b) + [K:\Q](h(g)+[K:\Q]\log(2)) = (\ell(K)+\ell_K(b))^{O(1)}.
$$
A sparsely given polynomial $f=\sum_{i=1}^t
a_ix^{\alpha_i}y^{\beta_i}\in K[x,y]$ is then encoded by the list of
its (non-zero) coefficients and corresponding exponents, and its
length relative to $K$ is
$$
\ell_K(f):=\sum_{i=1}^t (\ell_K(a_i) + \ell(\alpha_i) +
\ell(\beta_i)).
$$
Note that the input data is specified by $f$ {\it and} $K$, and so
the input length is $\ell(K) +\ell_K(f)$.
We have that
$$
t, \log_2(\deg f) \le \ell(f) \quad \mbox{ and } \quad
h(f) \le \ell_K(f) + [K:\Q](h(g)+[K:\Q] \log(2)) = (\ell(K)+\ell_K(f))^{O(1)}.
$$

\medskip
When the input of our algorithms comprises an inclusion
$K\hookrightarrow L$ of number fields,
$L$ is described as an extension of $K$ by a monic irreducible polynomial
$k(z) \in\OO_K[z]$ such that
 $L=K(\vartheta)$ for a root $\vartheta$ of $k$; this polynomial is represented in a dense way.
A polynomial $p\in L[x,y]$ in the output is
then encoded by the (dense) list of its coefficients
with respect to  the product basis $(\theta^j\vartheta^k)_{{0\le j\le \delta-1}, 0\le k\le \gamma-1}$ of $L$
over $\Q$; here we set $\gamma:=[L:K]$.
Note that for an element $b\in K$ in the base field encoded as
$b=b_0+ \cdots + b_{\delta-1}x^{\delta-1}$ with respect to the given
basis of $K$ over $\Q$,
its encoding with respect to the product base will be the same and so
$$
\ell_L(b) \le [L:K]\,\ell_K(b)
$$
since we have to count the zero coefficients corresponding to the monomials
$\theta^j\vartheta^k$ with $k\ge 1$.
In particular $\ell_L(f) \le [L:K]\,\ell_K(f)$ for $f\in K[x,y]$.

\medskip For the absolute factorization algorithm for $f\in K[x,y]$, the output irreducible polynomials
$p_i\in \Qbarra[x,y]$ are encoded by $(L_i,p_i)$, where $L_i$
consists in the minimal extension of $K$ such that $p_i\in L_i[x,y]$
(we observe that this encodes a full set $(\sigma(p_i))_{\sigma:
K\hookrightarrow \Qbarra}$ of $[L_i:K]$ conjugate factors of $f$).
The couple $(L_i,p_i)$ is encoded by a monic irreducible polynomial
$k_i(z) \in\OO_K[z]$ such that
 $L_i=K[z]/(k_i(z))$, and $p_i$ is given by its coefficients.

\subsection{Binomial factors} \label{Binomial factors}

The computation of the irreducible factors of a bivariate polynomial
that are binomials of more generally
products of binomials
 can be reduced to the univariate case
as we show in this section.
We first observe that an irreducible polynomial $p\in K[x,y]$ is a product of binomials
if it has one of the following forms:
\begin{equation}\label{binfactors}
p(x,y)=\prod_{\sigma} (x^\alpha - \sigma(\xi) y^\beta) \quad
\mbox{or} \quad p(x,y)=\prod_{\sigma} (x^\alpha y^\beta-
\sigma(\xi)),
\end{equation}
where $\alpha,\beta\ge 0$ are not $0$ simultaneously, $\xi\in \Qbarra$ and
where   $\sigma:K(\xi)\hookrightarrow \Qbarra$ runs over all $K$-embeddings
of $K(\xi)$ in $\Qbarra$.

We have the following results:

\begin{lem} \label{mult}
Let   $ \alpha, \beta, n\in \N$, $\xi\in \Qbarra^\times$ and $f\in
\Qbarra[x,y]$ be given. Set  $z$ for a new variable  and denote
 by $g\in \Qbarra[x,y,z]$  the remainder of the division  with respect to the variable $x$  of $f(x,y)$ by the monic
polynomial $x^\alpha-zy^\beta$.  Then
 $$(x^\alpha-\xi\,y^\beta)^n\div f(x,y) \quad \iff \quad
(z-\xi)^n \div g(x,y,z).$$
\end{lem}

\begin{proof}
Consider the ring
$$
 A:=\Qbarra[x,y^{\pm 1},z]/(x^\alpha-zy^\beta).
$$
We have that $x^\alpha-\xi y^\beta = (z-\xi)y^\beta$ in $A$, and,
since $y$ is invertible, we have the following equality of ideals
$$((x^\alpha-\xi y^\beta)^n) \ = \ ((z-\xi )^n)  \ \mbox{in }  A.$$
We call this ideal $I$. By definition $f=g$ in $A$ and so $f\in I$
if and only if $g\in I$, that is
 $$(x^\alpha-\xi\,y^\beta)^n\div f(x,y)  \quad \mbox{\rm in } A
\quad \iff \quad
(z-\xi)^n \div g(x,y,z) \quad \mbox{\rm in } A.$$
We have to show that we can take out the words ``in $A$'' from the above statement.

We observe that there is a natural identification
$A=\Qbarra[x,y^{\pm1}] $. Therefore, $$ (x^\alpha
-\xi\,y^\beta)^n\div f \mbox{ in } \ A  \ \iff \ (x^\alpha-\xi
\,y^\beta)^n\div f  \mbox{ in } \Qbarra[x,y^{\pm 1}]\ \iff \
(x^\alpha-\xi \,y^\beta)^n\div f\mbox{ in } \Qbarra[x,y]$$ since $y$
is prime to $x^\alpha -\xi \,y^\beta$.

We have a second identification
$$
A=\bigoplus_{j=0}^{\alpha-1} \Qbarra[y^{\pm 1},z] \cdot x^j,
$$
and therefore  $$(z-\xi)^n\div g  \mbox{ in } A \ \iff \
(z-\xi)^n\div g \mbox{ in } \Qbarra[x,y^{\pm 1},z]\ \iff \
(z-\xi)^n\div g \mbox{ in } \Qbarra[x,y,z]$$ since $y$ is prime to
$z-\xi$.
\end{proof}

\begin{cor} \label{cor_mult}
With the same notations than in the previous lemma, let $K$ be a number field
and suppose that $f\in K[x,y]$.
Set
$$
p(x,y):=\prod_{\sigma}
(x^\alpha - \sigma(\xi) y^\beta) \in K[x,y] \quad \mbox{and} \quad
q(z):=\prod_{\sigma}
(z - \sigma(\xi) ) \in K[z]
$$
where $\sigma$ runs over all $K$-embeddings of $K(\xi)$ in $\Qbarra$, then
 $$p(x,y)^n \div f(x,y) \quad \iff \quad
q(z)^n
 \div g(x,y,z).$$
\end{cor}

\begin{proof}
The polynomials $x^\alpha - \sigma(\xi) y^\beta$ for different $\sigma$'s are relatively prime,
and the same is true for the
polynomials $z - \sigma(\xi)$. Hence
$p(x,y)^n \div f(x,y)$ if and only if
$(x^\alpha - \sigma(\xi) y^\beta)^n\div f(x,y)$ for all $\sigma$ if and only if
$(z - \sigma(\xi))^n \div g(x,y,z)$ for all $\sigma$ if and only if
$q(z)^n\div g(x,y,z)$.
\end{proof}

The algorithm to compute the irreducible  factors  of $f\in K[x,y]$, of degree bounded by $d$,
that  are product of
binomials is now clear:

We are looking for factors $p(x,y)\in K[x,y]$ of degree $\le d$ of
one of the forms in (\ref{binfactors}). The cases $\xi=0$,
$\alpha=0$ or $\beta=0$ reduce directly to the univariate case where
we apply Lenstra's algorithm \cite[Thm]{Len99b} to the corresponding
content of $f$.

 So we can restrict ourselves to the cases when $\xi\in \Qbarra^\times$ and $\alpha,\beta\in \N$.
We consider first the factors of the first form in
(\ref{binfactors}).

We fix $1\le \alpha,\beta \le d$, and we set
$
g:=g_{\alpha,\beta}\in K[x,y,z]
$
for the remainder of dividing $f$ (with respect to $x$) by $x^\alpha-z y^\beta$
($g$ depends only on $f$ and $\alpha,\beta$).
It is easy to compute $g$ by Euclidean division:
$$
g(x,y,z) =\sum_{i=1}^t a_i x^{\alpha_i \bmod \alpha}(z\, y^\beta)^
{\lfloor \alpha_i/\alpha\rfloor} y^{\beta_i},
$$
so that $g$ is as sparse as $f$. We write
$$g(x,y,z)=\sum_{i,j}g_{i,j}(z)x^iy^j$$ and observe that an irreducible factor $q\in K[z]$ satisfies
$q^n\div g \iff q^n\div g_{i,j} \quad  \mbox{for all} \ i,j$, where
there are at most $t$ non-zero polynomials $g_{i,j}$, and each of
them is as sparse as $f$, with coefficients obtained as the sum of
at most $t$ coefficients of $f$.

We compute all irreducible factors $q\in K[z]$ of $g$ of degree
bounded by $ d/\max\{\alpha,\beta\}$ and  their corresponding
multiplicities, by examining the common irreducible factors (and
their multiplicities) of all the $g_{i,j}$'s. This is done again
applying Lenstra's univariate algorithm.

Since the irreducible polynomial $q$ is of the form $q=\prod_\sigma
(z-\sigma(\xi))$, the corresponding irreducible factor $p$ of $f$ is
then derived as
$$
p(x,y)=(y^\beta)^{\deg(q)} q(x^\alpha y^{-\beta}),
$$
where $\deg(p) = \max\{\alpha,\beta\} \cdot \deg(q) \le d$.

Corollary~\ref{cor_mult} certifies that for given $\alpha,\beta$, we
obtain in this way all irreducible factors of $f$ of degree $\le d$
of the first form in (\ref{binfactors}), as well as their
multiplicities.

\smallskip
For the factors in (\ref{binfactors}) of the second form,
we proceed similarly, by considering the remainder $g\in K[x,y^{\pm 1},z]$  of dividing $f$
(with respect to $x$) by
$x^\alpha y^\beta -z$.
We observe that the corresponding  extensions of Lemma~\ref{mult} and
Corollary~\ref{cor_mult}
hold.
In this case, $p$ is  derived from the factor $q\in K[z]$ of $g$ as
$
p(x,y)= q(x^\alpha y^\beta)
$.

\bigskip
The  algorithm described above yields the following result:

\begin{teo} \label{cyclotomic}
 There is a deterministic algorithm that, given
 $f\in K[x,y]$ and
$d\ge 1$, computes all irreducible factors of $f$ in $K[x,y]$ of
degree $\le d$ which are products of binomials, together with their
multiplicities, in $\big(d\cdot(\ell(K) +\ell_K(f))\big)^{O(1)}$ bit
operations.
\end{teo}

\begin{proof}
We have already established that the previous algorithm gives these
factors and their multiplicities. Its running time is estimated as
follows: for each pair $\alpha,\beta$, we are applying Lenstra's
algorithm $\le t$ times to the polynomials $g_{i,j}$ of sparse
length $\ell(g_{i,j})=O(\ell(f))$, in order to compute  their
irreducible  factors of degree $\le d/\max\{\alpha,\beta\}$ and
their multiplicities. This task is done in $\big(d\cdot(\ell(K)
+\ell_K(f))\big)^{O(1)}$ bit operations. Since there are at most
$d^2$ pairs $\alpha,\beta$, the total bit cost of the algorithm
remains of order $\big(d\cdot(\ell(K) +\ell_K(f))\big)^{O(1)}$.
\end{proof}

\subsection{Rational factorization} \label{Rational factorisation}

The search of all the low degree factors of a sparse $f\in K[x,y]$
is done by decomposing it as a sum of short pieces, as in the
previous papers \cite{CKS99,Len99b,KaKo05}. For given
$\Delta_x,\Delta_y\ge 0$, these pieces have to be separated by a
distance (``gap'') of at least $\Delta_x$ in the $x$-direction or
$\Delta_y$ in the $y$-direction. This is done here by decomposing
$f$ first  with respect to the $y$-exponents, then with respect to
the $x$-exponents.

\medskip
Let $f=\sum_{i=1}^t a_ix^{\alpha_i}y^{\beta_i}$ and suppose that
the monomials are already ordered so that
$\beta_1\le \beta_2\le  \cdots \le \beta_t$.
Then we determine
$$
\ell_0:=0 < \ell_1< \cdots < \ell_s <\ell_{s+1} =t
$$
subject to the conditions
$$\beta_{i+1}-\beta_{i} 
<\Delta_y \ \mbox{ for } \ell_{j}+1 \le i\le \ell_{j+1},  0\le j\le
s, \quad \mbox{ and } \quad
\beta_{\ell_j+1}-\beta_{\ell_j} 
\ge \Delta_y \  \mbox{ for } 1\le j\le s,$$
namely we split the $y$-exponents $ \beta_1,\dots,\beta_t$
into subsets so that consecutive exponents in the same subset are at distance $< \Delta_y$
and between different subsets there is a
gap of length $\ge \Delta_y$.
Set
$$
r_j:=\sum_{i=\ell_j+1}^{\ell_{j+1}}a_{i}x^{\alpha_{i}}y^{\beta_{i}-\beta_{\ell_j+1}}
\ \mbox{ for } 0\le j\le s
\quad \mbox{ so that }
\quad f=y^{\beta_{\ell_0+1}}r_0 +
 y^{\beta_{\ell_1+1}}r_1  + \cdots + y^{\beta_{\ell_{s}+1}}r_s.
$$
Next we do the same procedure over each $r_j$ with respect to $\Delta_x$: first we reorder the
monomials applying a permutation $\tau$ so that
$$
r_j= \sum_{i=\ell_j+1}^{\ell_{j+1}}a_{\tau(i)}x^{\alpha_{\tau(i)}}y^{\beta_{\tau(i)}-\beta_{\ell_j+1}}
$$
and
$
\alpha_{\tau(\ell_j+1)} \le \alpha_{\tau(\ell_j+2)} \le\cdots
\le \alpha_{\tau(\ell_{j+1})} $.
Then for each $0\le j\le s$
we sub-split this set of $\ell_{j+1}-\ell_j$ exponents
into subsets such that the consecutive $x$-exponents in the same subset are at distance $< \Delta_x$,
and between different subsets there is a
gap of length $\ge \Delta_x$.
Using this, we decompose $r_j$ into pieces
$$
r_j =  x^{\zeta_{0,j}}r_{0,j}
+ \cdots + x^{\zeta_{t_j,j}}  r_{t_j,j}
$$
for some exponents $\{\zeta_{i,j} : 0\le j\le s, 0 \le i\le t_j \}
\subset \{\alpha_1,\dots, \alpha_t\}$
that we do not explicit to avoid useless
proliferation of indexes.

\medskip
Each $r_{i,j}$ is (up to a monomial) some part of $r_j$, which in time
is (up to a monomial) some part of $f$.
We arrive in this way to a list
of $k\le t$ non-zero polynomials $f_1,\dots,f_k$
(after rewriting the $r_{i,j}$'s
into $f_i$'s) such that
\begin{equation}\label{subdivision}
f=x^{ \gamma_1}y^{\delta_1}f_1+x^{ \gamma_2}y^{\delta_2}f_2+\cdots
+ x^{\gamma_k}y^{\delta_k}f_k;
\end{equation}
and by construction for $ 1\le i\le k$,
$$
\ell_K(f_i)\le \ell_K(f), \quad \deg_x(f_i)<
(t-1)\,\Delta_x, \quad  \deg_y (f_i)< (t-1)\,\Delta_y
$$
and for $i\ne j$ we have that
\begin{align*}
\mbox{ either } & \quad \gamma_j -\gamma_i - \deg_x(f_i) \ge \Delta_x
\quad \mbox{ or }  \quad \gamma_i -\gamma_j - \deg_x(f_j) \ge \Delta_x \\
\mbox{ or }  & \quad \delta_j -\delta_i - \deg_y(f_i) \ge \Delta_y \quad \mbox{ or }
\quad \delta_i -\delta_j - \deg_x(f_j) \ge \Delta_y.
\end{align*}
We have decomposed $f$ in $\le t$ pieces of controlled degree and separated
by a gap of length $\ge \Delta_x$ in the $x$-direction or $\ge \Delta_y$ in the $y$-direction.

\bigskip
The computation of the irreducible factors
of $f$ of degree $\le d$ is then clear.
First we compute
a constant
$c$ such that $h_1(f)\le c$ in time
$(\ell(K)+ \ell_K(f))^{O(1)}$, as in~\cite[Prop.3.6]{Len99b}.
We set
$$
\Delta_x:=\Delta_y :=  \Delta=
5^6\cdot[K:\Q] \cdot d \cdot \log^3 (16[K:\Q] d)\cdot
c.
$$
Applying Corollary~\ref{gap_aritmetico} we infer that for $f=x^{
\gamma_1}y^{\delta_1}f_1+x^{ \gamma_2}y^{\delta_2}f_2+\cdots +
x^{\gamma_k}y^{\delta_k}f_k$ as in~(\ref{subdivision}), then for
$p\in K[x,y]$ that is not a  cyclotomic polynomial, we have
$$
p \div f \iff p\div f_i \quad \mbox{ for all } i.
$$
The procedure consists on computing first the cyclotomic factors
together with their multiplicity, by using the algorithm in
Subsection~\ref{Binomial factors}. For the other factors, we compute
them as the common factors of the $f_i$'s, by using any
polynomial-time algorithm for factoring dense bivariate polynomials
over a number field, for instance~\cite[Thm~3.26]{Len87}. Therefore
we obtain the following result:

\begin{teo} \label{numberfact}
There is a deterministic algorithm that, given
 $f\in K[x,y]$ and
$d\ge 1$, computes all irreducible factors of $f$ in $K[x,y]$ of
degree $\le d$, together with their multiplicities, in
\\$\big(d\cdot(\ell(K) +\ell_K(f))\big)^{O(1)}$ bit operations.
\end{teo}

\begin{proof}
We have already established that the previous algorithm gives all
these factors and their multiplicities. We estimate its running
time. We show that the degree of $f_i$ for all $i$, $1\le i\le k$,
in the decomposition (\ref{subdivision}) is polynomial in the input
size. This is a consequence of our estimate for the gap length:
$$
\ell(f_i)\le \ell(f) \ \  \mbox{ and } \ \ \deg_x (f_i),\deg_y
(f_i)<
(t-1)\,\Delta=
O(t \cdot([K:\Q] \cdot d)^{1+\varepsilon} \cdot c)
=
\big(d\cdot(\ell(K) +\ell_K(f))\big)^{O(1)}.
$$
Then we apply to each $f_i$ a polynomial-time algorithm for factoring dense
bivariate polynomials over $K$, which would do the task
in
$\big(d\cdot(\ell(K) +\ell_K(f))\big)^{O(1)}$ bit operations.
Since the number of $f_i$'s is  at most $t\le \ell$, the total complexity remains of the same order.
\end{proof}

If for an input polynomial $f\in K[x,y]$ we are interested  in its
 factors in an extension $L$,  we can compute them by just including $f$ into $L[x,y]$
 and then applying the above algorithm over $L$;
its cost would be of
$\big(d \cdot (\ell(K) +\ell_K(f)+ \ell_K(L)\big)^{O(1)}$ bit operations.

\smallskip
We note that here, for the factors which are products of binomials but not cyclotomic,
we have the choice of computing them either by
 reduction to the univariate sparse case of Theorem \ref{cyclotomic} or by reduction to the
 dense bivariate case.

\subsection{Absolute factorization} \label{Absolute factorisation}

Given a polynomial $f\in K[x,y]$, we can apply Corollary
\ref{gap_geometrico} to  extend the previous algorithm to the
computation of all irreducible factors of $f$ over $\Qbarra$, of
degree bounded by $d$, except the binomial ones. We assume that the
input $f$ is encoded in $K[x,y]$ and as before we compute a constant
$c$ such that $h_1(f)\le c$ in time $(\ell(K)+ \ell_K(f))^{O(1)}$,
then we set
$$
\Delta_x:=\Delta_y :=  \Delta=
2^{70} \cdot d \cdot \log^5(d+2) \cdot c.
$$
Corollary~\ref{gap_geometrico} implies that for the associated decomposition
$f=x^{ \gamma_1}y^{\delta_1}f_1+x^{ \gamma_2}y^{\delta_2}f_2+\cdots
+ x^{\gamma_k}y^{\delta_k}f_k$ as in (\ref{subdivision}),
 any  irreducible $p\in \overline\Q[x,y]$
that is not of the form
$$p(x,y)=x^\alpha - \xi y^\beta \ \mbox{ or } \
p(x,y)=x^\alpha y^\beta - \xi,
$$
satisfies
$$
p \div f \iff p\div f_i \quad \mbox{ for all } i.
$$
Now we need to determine the common factors of  the $f_i$'s  over
$\overline\Q[x,y]$ and  their multiplicity. In order to do this,
 we first factor completely each of the $f_i$ over $K[x,y]$ by
applying any dense polynomial-time bivariate factorization algorithm
over $K$. An irreducible factor $p\in \overline\Q[x,y]$ of $f$ will
necessarily divide a common irreducible factor $q\in K[x,y]$ of all
the $f_i$'s. Thus it is enough to keep all common
 irreducible factors  $q\in K[x,y]$ of all the $f_i$'s and
their multiplicities, and then  to  factor them in
$\overline\Q[x,y]$ by applying  any polynomial-time algorithm for
factoring dense bivariate polynomials over $\overline\Q$, for
instance~\cite[Theorem~11]{Kal95}. We only keep those factors in the
output which are of degree $\le d$ and which are not binomials. We
proceed in this way in order to avoid comparing irreducible factors
in $\Qbarra[x,y]$ of different $f_i$'s, that can, although equal, be
described in different field extensions.

\begin{teo} \label{absolute}
There is a deterministic algorithm that, given $f\in K[x,y]$ and $d
\ge 1$, computes all irreducible factors of $f$ in $\Qbarra[x,y]$ of
degree $\le d$, together with their multiplicities, except for the
binomial ones, in $\big(d\cdot(\ell(K) +\ell_K(f))\big)^{O(1)}$ bit
operations.
\end{teo}

\begin{proof}
As with the previous one, the complexity of this algorithm is
estimated in  $\big(d \cdot(\ell(K) +\ell_K(f))\big)^{O(1)}$ bit
operations, because we have to factor $\le t$ polynomials $f_i$ of
degree polynomially bounded in the input length to find all possible
$q$, which are of input length $\ell_K(q)=\big(d \cdot(\ell(K)
+\ell_K(f))\big)^{O(1)}$ and at most the same quantity, and then to
factor them in $\Qbarra[x,y]$.
\end{proof}

\subsection{A practical improvement: adaptive gap methods} \label{adaptive}

The practical efficiency of the proposed algorithms depends
essentially on the length $\Delta$ defining the gap in $f$: the
degree of the pieces $f_i$  depends  on $\Delta$, and if this degree
is large, the dense factorization algorithm will be clearly slower.
In other words, the smaller the gap length $\Delta$ is, the faster
the algorithm  works. Since the gap is  proportional to the inverse
of the essential minimum, the greatest the essential minimum, the
faster the algorithm.

\medskip
There are some special situations where we can get better bounds, for instance
for linear factors $p(x,y)=ax+by+c$ with integer coefficients, as in~\cite{KaKo05}.

The Mahler measure of a polynomial is bounded from below by the
Mahler measure of any of its facet polynomials. Hence for $a,b,c\in
\Z$ relatively prime numbers such that $a\cdot b\cdot c\ne 0$, we
have that
$$
m(ax+by+c)\ge \max\{ m(ax+by) ,m(by+c),m(ax+c)\} = \log \max\{|a|,|b|,|c|\}
$$
as it can be proved that the Mahler measure of a binomial coincides with its height.
The theorem of successive minima then implies
$$
\mu^\ess(Z(ax+by+c)) \ge \frac12 \log \max\{|a|,|b|,|c|\}= \frac12 h(p).
$$
The only case for which this lower bound is meaningless is when
$a,b,c=0,\pm 1$. (When $a$, $b$ or $c$ vanish, we reduce easily to
the univariate case so we do not consider it here.) When $a,b,c=\pm
1$, Zagier's theorem~\cite{Zag93}, see also Subsection \ref{plane
curves}, shows that $ h(\xi) \ge 0.1911$. Hence
$$
\mu^\ess(Z(ax+by+c)) \ge \left\{
\begin{array}{ll}
 \log(\xi_0) = 0.1911 &\quad \mbox{ if } a,b,c=\pm 1 \\[2mm]
 h(p) \ge \frac{ \log(2)}{2} = 0.3465 & \quad \mbox{ otherwise. }
\end{array}
\right.
$$
which improves the bound $\log(1.045)\approx  0.0440$ proposed in \cite{KaKo05}.

Note that in this case the gap size associated with $p=ax+by+c$ gets smaller
as the coefficients of $p$ tend to infinity.
Therefore, a good strategy to make the algorithm more efficient   might be to exclude a finite
number of candidates by testing them as factors of $f$ (using a rough estimate for their gap length),
 and then use a much smaller gap length to find the rest of the factors by reduction to the dense case.

\end{document}